\documentclass[preprint,12pt, a4paper]{elsarticle}

\usepackage{amssymb,amsmath}
\usepackage{hyperref}
\usepackage{bm}
\usepackage{color}
\usepackage{booktabs}
\usepackage{wrapfig}
\usepackage{dirtree}
\usepackage{graphicx}
\usepackage{makecell}

\definecolor{green}{rgb}{0.0, 0.5, 0.0}

\journal{SoftwareX}

\begin{document}

\renewcommand{\labelenumii}{\arabic{enumi}.\arabic{enumii}}

\begin{frontmatter}

\title{rDSM---A robust Downhill Simplex Method software package for optimization problems in high dimensions}

\author[label1]{Tianyu Wang}
\author[label1]{Xiaozhou He\corref{mycorrespondingauthor}}\ead{hexiaozhou@hit.edu.cn}
\author[label2,label3,label1]{Bernd R.~Noack\corref{mycorrespondingauthor}}\ead{bernd.noack@szu.edu.cn}
\cortext[mycorrespondingauthor]{Corresponding author}

\address[label1]{School of Robotics and Advanced Manufacture, Harbin Institute of Technology, Shenzhen, 518055, Shenzhen,  P.~R.~China}

\address[label2]{College of Mechatronics 
and Control Engineering, Shenzhen University, Canghai campus,
518060 Shenzhen, P.~R.~China}

\address[label3]{Guangdong Province VTOL Aircraft Manufacturing Innovation Center,
518060 Shenzhen, P.~R.~China}

%

\begin{abstract}

The Downhill Simplex Method (DSM) is a fast-converging derivative-free optimization technique for nonlinear systems. 
However, the optimization process is often subject to premature convergence due to degenerated simplices or noise-induced spurious minima.  
This study introduces a software package for the robust Downhill Simplex Method (rDSM), which incorporates two key enhancements.
First, simplex degeneracy is detected and corrected by volume maximization under constraints. 
Second, the real objective value of noisy problems is estimated by reevaluating the long-standing points.
Thus, rDSM improves the convergence of DSM, and may increase the applicability of DSM to higher dimensions, even in the presence of noise.
The rDSM software package thus provides a robust and efficient solution for both analytical and experimental optimization scenarios.
This methodological advancement extends the applicability of simplex-based optimization to complex experimental systems where gradient information remains inaccessible and measurement noise proves non-negligible.

\end{abstract}


\end{frontmatter}

\begin{table}[!h]
\caption{Code metadata} 
\begin{tabular}{ l p{6.5cm} p{6.5cm} }
\toprule
\textbf{Nr.} & \textbf{Code metadata description} & 
\\
\midrule
C1 & Current code version & v1.0 \\
C2 & Permanent link to code/repository used for this code version & \url{https://github.com/tianyubobo/rDSM} \\
C3  & Permanent link to Reproducible Capsule & 
\url{https://github.com/tianyubobo/rDSM}
\\
C4 & Legal Code License   & CC-BY-SA
\\
C5 & Code versioning system used & git
\\
C6 & Software code languages, tools, and services used & MATLAB
\\
C7 & Compilation requirements, operating environments \& dependencies & 
 Microsoft Windows  \\
C8 & If available Link to developer documentation/manual & \url{https://github.com/tianyubobo/rDSM} \\
C9 & Support email for questions & wangtianyu@stu.hit.edu.cn \\
\bottomrule
\end{tabular}
\label{codeMetadata} 
\end{table}

\section{Motivation and significance}
The Downhill Simplex Method (DSM), originally formulated by Nelder and Mead in 1965~\cite{Nelder1965}, has served as a fundamental derivative-free optimization technique for multidimensional unconstrained optimization. 
The efficiency of the algorithm is attributable to its unique ability to handle non-differentiable objective functions.
This characteristic is especially beneficial in engineering applications where gradient-based optimization methods are not applicable.
Examples of applications include wind turbine problems~\cite{Xu2020}, structural engineering problems~\cite{SuHL2020,Wang2019CCDC}, civil engineering~\cite{Khalid2020}, material design engineering~\cite{HARDT2021,KOCI2016}, to cite a few examples.
However, the method has presented limitations in identifying optimal solutions when dealing with a large number of parameters and high-dimensional design spaces.

Numerous methods have been proposed in the literature to improve the performance of DSM in finding better optima. 
They can be broadly categorized into two groups: improving the internal mechanisms of DSM itself and combining with other optimization algorithms, i.e., hybrid methods. 
The initialization method and the selection of parameters constitute the primary focus of enhancing the meta-parameter tuning method of DSM.

Huang et al.~\cite{HUANG1998} applied a multi-start downhill simplex method for spatio-temporal source localization in magnetoencephalography, with the number of initializations typically ranging from 100 to 5000 automatically.
Kelley~\cite{Kelley1999} proposed a re-initialization methodology for the simplex as an alternative shrink step to prevent DSM from stopping at local extrema.
Koshel~\cite{koshel2002,koshel2005} demonstrated that optimizing the reflection, contraction, and expansion coefficients could reduce the number of iterations by up to 20\% to achieve the same performance. 
However, the effectiveness of such an approach in higher-dimensional problems remains untested. 
Gao and Han~\cite{Gao2012} provided a comprehensive analysis and recommendations concerning the selection of these coefficients in varying dimensions, validating their methodology in test functions with dimensionality greater than 100.

\begin{table}[h]
\centering
\caption{Table of past contributions to DSM.}\label{table: literature review}
\begin{tabular}{ccc}
\toprule
Notation      & Correct degenerated simplex & Noisy dataset \\ \midrule
Huang et al.~\cite{HUANG1998}  &  -   &  $\checkmark$     
\\
Kelley~\cite{Kelley1999}  &  -   &  -                           \\
Koshel~\cite{koshel2002,koshel2005}   &  -   &  -                           \\
Gao and Han~\cite{Gao2012}&  -   &  -                       \\
Luersen and  Le Riche~\cite{luersen2004}  &     $\checkmark$       &    -                     \\
Present study & $\checkmark$   &        $\checkmark$                 \\ \bottomrule
\end{tabular}
\end{table}

The integration of the DSM with other optimization algorithms to enhance the identification of better optima has also been the focus of other studies.
Pan and Wu~\cite{Pan1998} enhanced the robustness of DSM by incorporating simulated annealing strategies, achieving significant improvements in avoiding local extrema and reducing computational time. 
Similarly, Bangert~\cite{bangert2005} utilized DSM to optimize the meta-parameters of the simulated annealing algorithm, resulting in a 26.1\% improvement in the objective function compared to the results without DSM.  
Maehara and Shimoda~\cite{MAEHARA2013} proposed a hybrid methodology combining DSM with a genetic algorithm (GA), leveraging the strengths of both methods to compensate for their respective weaknesses.
Cornejo Maceda et al.~\cite{CornejoMaceda2021jfm,CornejoMaceda2023jfm} combine DSM with genetic programming to accelerate the learning of control laws for fluid flows in numerical simulations and experiments.
Furthermore, Wang Xin et al.~\cite{WangXin2023} applied DSM to minimize drag in fluidic pinball problems, using clustering techniques to control the corresponding drag.
Comparison of past methods is summarized in Table~\ref {table: literature review}.

Despite these efforts, a considerable gap remains in addressing both the convergence issues and the challenges posed by degenerated simplices. 
To bridge this gap, we add two targeted improvements to DSM, which is called the robust Downhill Simplex Method (rDSM):
    
\begin{itemize}
        \item \textbf{Degeneracy correction}. Degenerated simplices, where the vertices of the simplex become collinear or coplanar, compromise algorithmic efficiency and performance. This step rectifies dimensionality loss by restoring a degenerated simplex with $n-1$ or fewer dimensions to an $ n$-dimensional one, thereby preserving the geometric integrity of the search process.
        \item\textbf{Reevaluation.} This improvement prevents the simplex from getting stuck in a noise-induced spurious by reevaluating the cost function of the best point.
\end{itemize}

These modifications improve the convergence of DSM and may be beneficial for high-dimensional search spaces.
The two improvements are expected to extend the utility of DSM to complex optimization problems, mitigating prior limitations concerning degeneracy and noise.

The software architecture and functionalities are described in Sec.~\ref{Sec: description}. Sec.~\ref{Sec: example} gives the illustrative example and comparison of rDSM and DSM. The impact of this software package is described in Sec.~\ref{Sec: impact}; it is expected to benefit high-dimensional experimental optimization problems. Conclusions and outlook are provided in Sec.~\ref{Sec: conclusion}.

\section{Software description}\label{Sec: description}

The rDSM algorithm is a significant enhancement of the classic DSM that has been shown to improve convergence robustness and search capability in complex optimization landscapes.
This is achieved through the implementation of degeneracy correction and reevaluation strategies. 
The developed software package is implemented in MATLAB version 2021b.
It facilitates optimal search within high-dimensional spaces.
Given an initial point and objective function, the rDSM software package automatically gives the optimum, the learning curve and exports the information on the learning process (points explored, corresponding cost function, operations performed, and the counter values).

\subsection{Software architecture}

\begin{figure}[!htbp]
\includegraphics[width=0.98\textwidth]{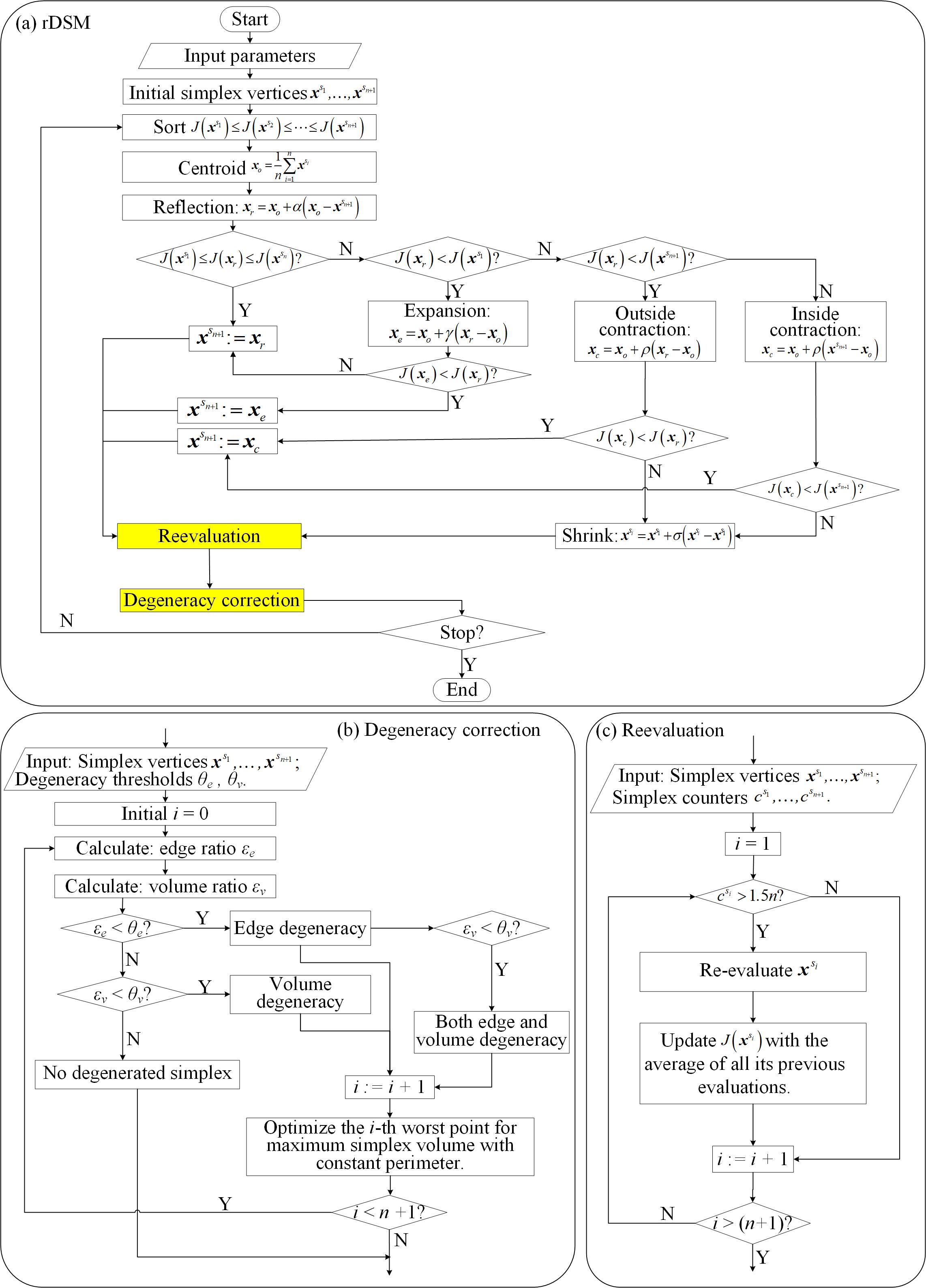}
\caption{Flowchart of rDSM. (a) shows the overview of rDSM, (b) and (c) give the detailed procedure of degeneracy correction and reevaluation, respectively.  }
\label{fig:rDSM}
\end{figure}

The rDSM improves the classic DSM with degeneracy correction and reevaluation as two enhancements, avoiding being trapped in a non-stationary point.
Fig.~\ref{fig:rDSM} shows the flowchart of the robust Downhill Simplex Method (rDSM) algorithm, demonstrating the framework design of the proposed optimization algorithm. It adds two steps of improvements for rDSM to the classic DSM framework design. 
The procedure of two steps, i.e., degeneracy correction and reevaluation, is shown in Fig.~\ref{fig:rDSM} (b) and (c), respectively.
The rDSM software package starts with the classic DSM procedure. 
For each iteration, the reflection, expansion, contraction, and shrink operations will be carried out based on the procedure shown in Fig.~\ref{fig:rDSM}. 
The two improvements are added subsequently. Mathematical and implementation details are given in Sec.~\ref{Sec: func}.
The rDSM can continue the optimization process by correcting the degenerated simplex with one or a few function iterations, thereby extending the exploration of the research domain. Reevaluation replaces the objective value of a persistent vertex with the mean of its historical costs, thereby enhancing the accuracy of the optimization.
Symbols describing the simplex and used in degeneracy correction and reevaluation are listed in Table~\ref{table: simplex quantity}.

\begin{table}[h]
\centering 
\caption{Table of simplex quantities and their notations.}\label{table: simplex quantity}
\begin{tabular}{cc}
\toprule
Notation         & Simplex quantities                                                                 \\ \midrule
$c^{s_i}$     &  Point $\bm{x}^{s_{i}}$ counter\\
$\bm{e}$         & Edge matrix  \\
$\bm{e}^i$       & $i$-th edge length                                                     \\
$J$              & objective function\\                                                     
$n$              & Searching space dimension               \\  
$P$              & Simplex perimeter 
\\
$V$              & Simplex volume
\\
$\bm{x}^{s_{i}}$ & $i$-th simplex point                              \\
$\bm{y}^{s_{n+1}}$         & Corrected point $\bm{x}^{s_{n+1}}$  
\\            
 \bottomrule
\end{tabular}
\end{table}

\begin{figure}[h]
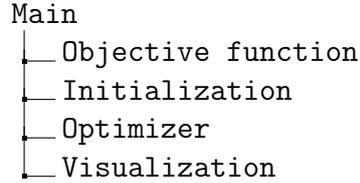

  \centering
    \begin{minipage}[h]{0.8\linewidth}
      \dirtree{%
          .1 Main\ .
            .2 Objective function .
            .2 Initialization . 
            .2 Optimizer .
            .2 Visualization . 
          }
  \end{minipage}
  \caption{Module structure of rDSM software package.\label{fig:structure}}
\end{figure}

The software package for rDSM is divided into four modules: (1) objective function, (2) initialization, (3) optimizer, and (4) visualization. 
There is one folder with the same name as each module containing the corresponding scripts.
The structure of four modules is illustrated in Fig.~\ref{fig:structure}.
The ``objective function" module defines an objective function, which may call an external computational fluid dynamics (CFD) solver or run an experiment, acting as an interface that other modules will call. 
The ``initialization"  module generates the initial simplex and operation parameters for optimization.
The ``optimizer"  module presents the iteration procedure of the optimization.
Finally, the ``visualization" module gives a figure to show the simplex iteration history and the learning curve.

\subsection{Objective function}
\begin{sloppypar}
The ``objection function" module defines the function to minimize.
One example of an objective function is given in the `\texttt{/ObjectiveFunction/test\_function.m}' file. 
It is a two-dimensional linear gradient function with an obstacle. 
Users are invited to implement their objective function. 
\end{sloppypar}
\subsection{Initialization}

\begin{sloppypar}
The ``initialization" module generates the initial simplex and initial operation parameters.
The default value for the initial coefficient of the first simplex is 0.05, which can be set a little larger for higher-dimensional problems.
The operation parameters are initialized in the file  `\texttt{/Initialization/DSM\_parameters\_N().m}'. 
The default reflection, expansion, contraction, and shrink coefficients are 1, 2, 0.5, and 0.5, respectively. 
In the rDSM software package, we introduce two parameters, the edge and volume thresholds, as a criterion to start correcting the degenerated simplex, details are in Sec.~\ref{sec:degeneracy correction}.
Users can change these parameters based on their optimization problem.
As suggested in~\cite{Gao2012}, the reflection, expansion, contraction, and shrink coefficients could be a function of the dimension $n$ of the search space, especially for $n>10$. 
All the notations and default values of the aforementioned parameters are summarized in Table~\ref{table: paramters}.
\end{sloppypar}
\begin{table}[h]
\centering
\caption{Table of rDSM parameters, notations, and default values.}\label{table: paramters}
\begin{tabular}{ccc}
\toprule
Parameter         & Notation   & Default value \\ \midrule
Reflection coefficient  & $\alpha$   & 1             \\
Expansion coefficient   & $\gamma$    & 2             \\
Contraction coefficient & $\rho$     & 0.5           \\
Shrink coefficient      & $\sigma$   & 0.5           \\
Edge threshold     & $\theta_e$ & 0.1           \\
Volume threshold   & $\theta_v$ & 0.1           \\ \bottomrule
\end{tabular}
\end{table}

\subsection{Optimizer}\label{Sec: func}
\begin{sloppypar}
The ``optimizer" module is the basis of this software package. 
All the function scripts for DSM and rDSM optimization are given in the folder `\texttt{/Optimizer}'.
The `\texttt{/Optimizer/DSM.m}' script is a re-implementation of the \texttt{fminsearch} function in MATLAB.
\texttt{/Optimizer/rDSM.m} is the optimizer introduced in this study. 
\end{sloppypar}

\subsubsection{Degeneracy detection}

The degeneracy correction step is used to correct the degenerated simplex.
An $n$-dimensional simplex is defined as the convex hull spanned by $n+1$ affinely independent points.
It is said to be \textit{degenerated} if the simplex is stretched such that one or more of its dimensions are negligible compared to the others.
This may happen, for example, when the simplex goes through a valley and the vertices of the simplex get aligned.
In this scenario, the degenerated simplex will be corrected by moving the worst point based on its objective value to maximize the volume and preserve its perimeter.

\begin{figure}[h]
\includegraphics[width=0.9\textwidth]{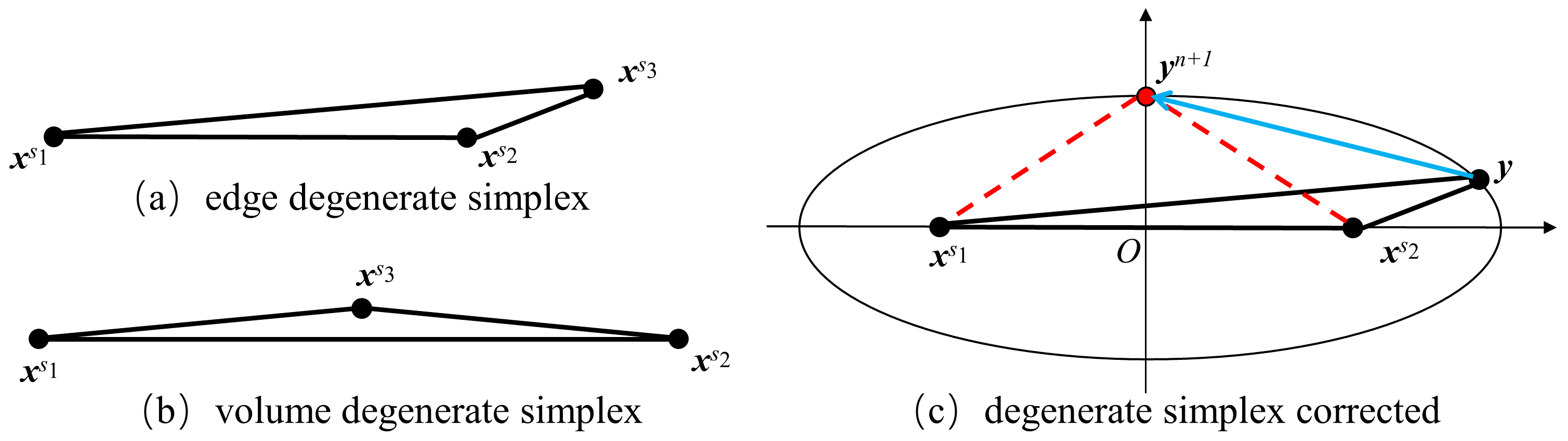}
\caption{
The two degenerated simplex types: (a) edge-degenerated simplex, (b) volume-degenerated simplex. and (c) the corrected degenerated simplex. The black solid lines are the edges of the degenerated simplex, and the red dashed lines present the corrected edges. The blue line shows the moving direction of the vertex to be corrected. $\bm{y}^{s_3}$ is the corrected point of $\bm{x}^{s_3}$.}

\label{fig:degenerate simplex}
\end{figure}

Let's consider a two-dimensional (2D) problem as an illustrative example.
It can be observed that when going through a valley, the operations of DSM tend to flatten the simplex, transforming it into a shape closer to a line segment, i.e., a 1-D simplex. 
This phenomenon can manifest in two distinct types as documented in ~\cite{luersen2004} and visualized in Fig.~\ref{fig:degenerate simplex}, (a) the edge degenerated simplex and (b) the volume-degenerated simplex. 
Because of the unique characteristics associated with each type of degeneracy, it is essential to utilize two separate criteria for their accurate identification.
The first case shown in Fig.~\ref{fig:degenerate simplex} (a) can be identified by comparing the shortest edge of the simplex with its largest edge.
Precisely, we compute the ratio $\varepsilon_e$ between the minimum edge length and the maximum edge length of this simplex.
If $\varepsilon_e$ is smaller than the edge threshold $\theta_e$, i.e.,
\begin{equation}\label{eq:dege}
    \varepsilon_e = \frac{\min_{i=1,n} \left \| \bm{e}^i \right \| }{\max_{i=1,n} \left \| \bm{e}^i \right \| } <\theta_e,
\end{equation}
the simplex is said to be \textit{edge-degenerated}. 
Symbol $\bm{e}^i$ represents the $i$-th edge vector, $\left \| \cdot  \right \| $ is the Euclidean norm, so that $\left \| \bm{e}^i  \right \|$ is the length of this edge vector. 
The default value for $\theta_e$ is 0.1, which means the largest edge of the simplex is 10 times larger than the shortest edge.
However, this criterion can not identify the degenerated simplex shown in Fig.~\ref{fig:degenerate simplex} (b), whose $\varepsilon_e$ is close to 0.5.
Therefore, a second criterion is given comparing the volume of the simplex and its edge length.
In practice, if the ratio $\varepsilon_v$ is smaller than thresholds $\theta_v$, i.e.,
\begin{equation}\label{eq:vol}
  \varepsilon_v = \sqrt[n]{\frac{\det\left [ \bm{e} \right ] }{\prod_{i}\left \| \bm{e}^i \right \| }} <\theta_v,
\end{equation}
the simplex is called \textit{volume-degenerated}, which means that, compared to the edge length, the volume of this simplex is too small. 
In Eq.~\eqref{eq:vol}, $\det\left [ \bm{e} \right ]$ denotes to the volume of the simplex, and symbol $\bm{e}$ represents the edge matrix of this simplex.
The default value for $\theta_v$ is 0.1.
By Eq.~\eqref{eq:vol}, we cannot identify the case in Fig.~\ref{fig:degenerate simplex} (a).
One example in 2D is that a simplex established by points (0,0), (1000,0), and (0,1), which is not \textit{volume-degenerated} but \textit{edge-degenerated}.
Therefore, both Eq.~\eqref{eq:dege} and Eq.~\eqref{eq:vol} are necessary for determining a degenerated simplex.
Both $\theta_e$ and $\theta_v$ are user-adjustable and can be modified according to the dimensionality and complexity of the specific problem.



\subsubsection{Degeneracy correction}\label{sec:degeneracy correction}
To correct the degenerated simplex, we fix its perimeter and maximize its volume.
In practice, we move the vertex that performs the worst to maximize the volume with a constant perimeter constraint.
If there are several possible positions for the worst vertex, we select the closest point as the new vertex of the corrected simplex.
An example in 2D is shown in Fig.~\ref{fig:degenerate simplex} (c), the black solid lines construct a degenerated simplex, and the red dashed lines are the corrected edges.
The solid blue line shows the moving direction of the corrected simplex vertex.
The red dot $\bm{y}^{s_{3}}$ is the corrected vertex of the black dot $\bm{x}^{s_3}$.
In this 2D case, the corrected simplex vertex is on the ellipse whose foci are the two vertices $\bm{x}^{s_1}$ and $\bm{x}^{s_2}$.

The perimeter $P(\bm{x}^{s_1},\ldots,\bm{x}^{s_{n+1}})$ of a simplex is the sum of all edge lengths. 
The edge length is the Euclidean distance between the two points:
\begin{equation}
    P(\bm{x}^{s_1},\ldots,\bm{x}^{s_{n+1}}) = \sum_{i=1}^{n}\sum_{j=i+1}^{n+1}  \left \| \bm{x}^{s_i}-\bm{x}^{s_j} \right \|,
\end{equation}
where $ \left \| \cdot \right \| $ represents the Euclidean distance, $\bm{x}^{s_{i}},i=1,\cdots,n+1$ are the points constructing the simplex.

The volume of this simplex $V(\bm{x}^{s_1},\ldots,\bm{x}^{s_{n+1}})$ is defined in~\cite{CHO1991}:
\begin{equation}
    V(\bm{x}^{s_1},\ldots,\bm{x}^{s_{n+1}})=\frac{1}{n!} \left |\det\begin{bmatrix}
  \bm{x}^{s_1},&\ \dots ,  & \bm{x}^{s_{n+1}}\\
  1,& \dots , &1
\end{bmatrix}\right |.
\end{equation}
Therefore, the degeneracy correction of the simplex 
$\{\bm{x}^{s_1},\ldots,\bm{x}^{s_n}, \bm{x}^{s_{n+1}}\}$ is achieved by solving the following constrained maximization problem:
\begin{equation}\label{eq:dege co}
\left\{\begin{matrix}{\bm{y}}^{s_{n+1}} = \mathrm{arg} \; \underset{\bm{y}}{\max} \; V(\bm{x}^{s_1},\ldots,\bm{x}^{s_{n}},\bm{y}) \\
s.t.\\  P(\bm{x}^{s_{1}},\ldots,\bm{x}^{s_n},\bm{y}^{s_{n+1}}) = P(\bm{x}^{s_1},\ldots,\bm{x}^{s_n},\bm{x}^{s_{n+1}}),
\end{matrix}\right.
\end{equation}
where $\bm{y}^{s_{n+1}}$ is the optimized $n+1$ vertex.
In practice, the Eq.~\eqref{eq:dege co} is solved by the Newton–Raphson method.
It is worth mentioning that if the simplex is still degenerated after updating the worst point $\bm{x}^{s_{n+1}}$, the next worst point is optimized.
This process is repeated following the reverse objective value order until the simplex is no longer degenerated or all the points are optimized.

\subsubsection{Reevaluation}

The developed ``reevaluation module" aims to improve the convergence robustness of the classic DSM in the presence of noise.
Repeatedly misvaluing the objective value of a vertex can lead to a sub-optimum due to a misdirection of the optimization process.
To mitigate this, we introduce a counter $c^{s_i}$ for vertex $\bm{x}^{s_i}$ to track the number of DSM iterations that $\bm{x}^{s_i}$ remains in the processing simplex.
If $c^{s_i} \geq 1.5n$ for $\bm{x}^{s_i}$, we will \textit{reevaluate} $\bm{x}^{s_i}$, and its objective value will be replaced by the average of all the previous evaluations.
Taking the average of repeated evaluations reduces the effect of random errors, bringing the result closer to the true value.
Most vertices leave after about $n$ times. Vertex remaining more than $n$ times indicates the optimizer searches a limited domain around this vertex.
This slightly higher coefficient (1.5) balances the efficiency and accuracy. A larger one needs more computation resources, a smaller one makes almost every vertex need to be \textit{reevaluate}. Both make the optimizer rDSM less efficient and correct.
This reevaluation strategy facilitates continued optimization and enhances the likelihood of converging on the optimum.

\subsection{Visualization}

In the script \texttt{main.m}, we provide a visualization of the simplex iteration history for two-dimensional problems and a learning curve after optimizing.
The scripts are given in the folder `/Visualization'.

\begin{sloppypar}
The optimization process automatically generates three structured data archives (`SimplexHistory.txt', `PointsDatabase.txt', and `ReevaluationHistory.txt') in the folder `/Output'.
ASCII format files `SimplexHistory.dat' and `PointsDatabase.dat' are also provided.
The `SimplexHistory.txt' and `SimplexHistory.dat' files contain structured records of simplex vertices, including the simplex vertex index, 
the simplex indices, the operation of the corresponding simplex, and the number of times the vertices remained active during the optimization procedure. 
Similarly, `PointsDatabase.txt' and `PointsDatabase.dat' provide comprehensive data on vertex values, simplex vertex index, corresponding objective values, the simplex index, and the operations performed on the corresponding simplex.
Additionally, `ReevaluationHistory.txt' offers detailed information on parameter and objective values before and after reevaluation, ensuring traceability of optimization iterations. These files collectively serve as essential resources for reconstructing optimization trajectories and validating computational results.
\end{sloppypar}


\section{Illustrative examples}\label{Sec: example}

This section employs analytical functions to illustrate the functionality of rDSM. Firstly, the comparison results of minimizing the two-dimensional function by DSM and rDSM will be presented. Furthermore, a higher-dimensional test is provided in order to illustrate the performance of rDSM in high-dimensional space.

\subsection{Two-dimensional test functions}
In this section, we use a two-dimensional analytical function, i.e., a 2D linear gradient function with and without an obstacle as an objective function to minimize. 
The analytical objective function~\cite{Abraham2009}:
\begin{equation}\label{eq:linear grad}
    J (x_1,x_2) = -\frac{(x_1-x_2)}{4}+0.5,
\end{equation}
is defined on $\mathrm{\Omega} = \{(x_1,x_2) \in [-1,1]\times[-1,1]\}$.
On $\Omega$, $J$
is characterized by a global minimum at (1, -1) with corresponding objective value $J = 0$.
The objective value of obstacle is set as $10^3$ at $-1<x_1<0$ and $-1<x_2<0$.
During the DSM iteration process, points are assigned infinite objective values if they exceed the search domain due to operations of reflection, expansion, etc.

The performance of DSM and rDSM optimizers is evaluated by minimizing Eq.~\eqref{eq:linear grad} with three different cases: (a) DSM without obstacle; (b) DSM and (c) rDSM with obstacle.
The starting point for all cases is ($- 0.75$, $0.35$).
The iterative process is limited to 50 iterations, with a maximum of 100 evaluations. 
Both edge and volume degeneracy thresholds, i.e., $\theta_e$ and $\theta_v$, are determined based on prior tests on 2D problems,
considering the number of dimensions, the complexity of the optimization problem, etc..
Typically, $\theta_e$ and $\theta_v$ are set as small values suggested in~\cite{luersen2004}.
For this investigation, these parameters are assigned a value of 0.1.

\begin{figure}[!htbp]
\centering
\includegraphics[width=0.95\textwidth]{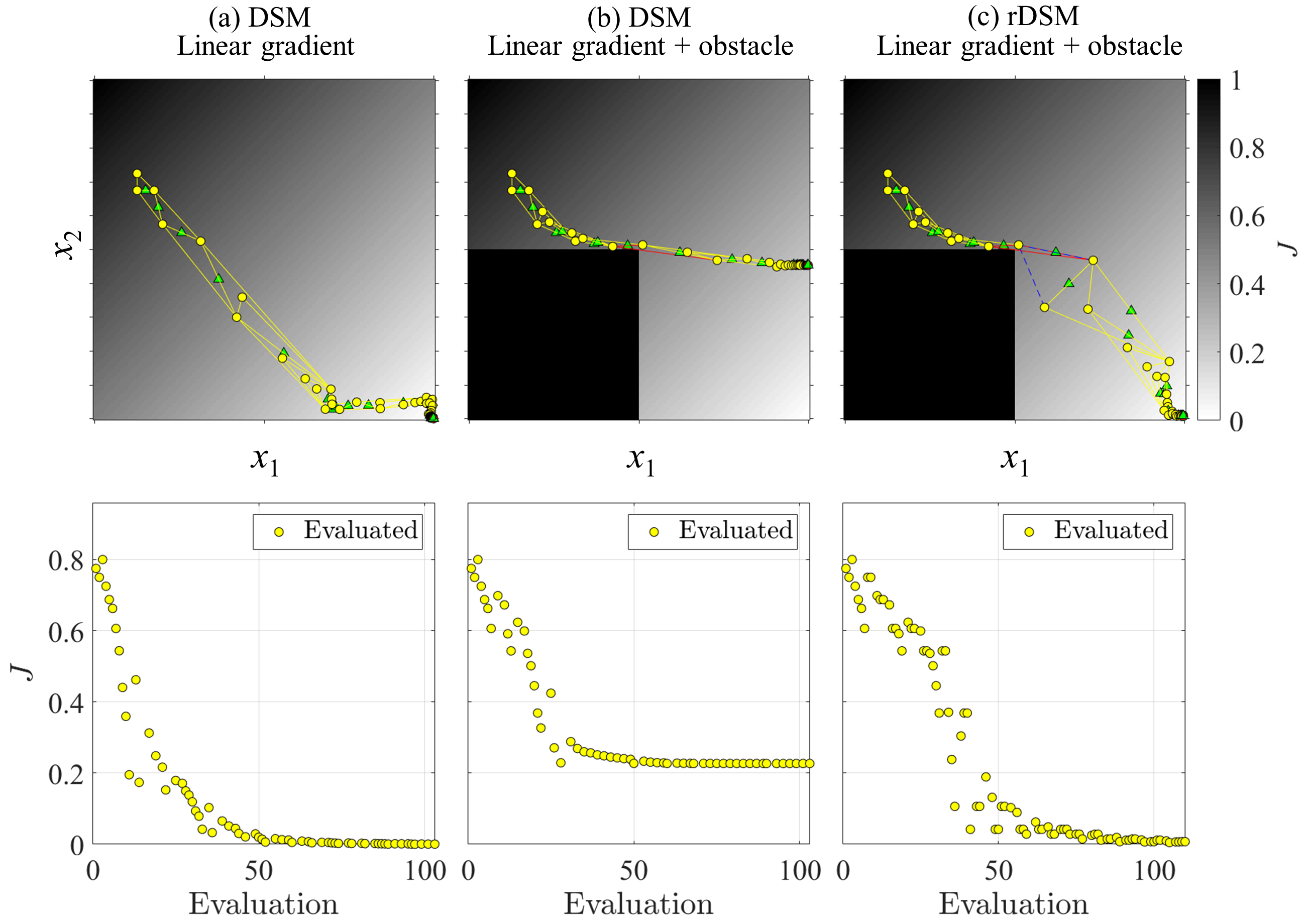}
\caption{Comparing the optimization process and corresponding learning curve of linear gradient with and without obstacle by DSM and rDSM.
}
\label{fig:degeneracy}
\end{figure}
The optimization results are shown in Fig.~\ref{fig:degeneracy}.
The first row shows the optimization process using the objective function value as a background for the three aforementioned cases. 
Simplex vertices throughout the optimization are denoted as yellow dots.
The centroids of each simplex, marked as green triangles, are included to show the learning process. 
Solid yellow lines present the simplices generated by DSM.
Solid red lines highlight the edges of degenerated simplices.
Blue dashed lines represent the corrected ones.
The second row displays the learning curve of corresponding cases, with evaluated points marked as yellow dots.

In the first column of Fig.~\ref{fig:degeneracy}, the DSM optimizer quickly converges at point (0.9999, -0.9998) with an objective value of $9.2146\times10^{-5}$.
The second column of Fig.~\ref{fig:degeneracy} illustrates the degenerated simplex during the iteration process and fails to reach the global minimum.
The endpoint for case (b) is (1, -0.0925) with a corresponding objective of 0.2269.
The simplex becomes flattened due to the obstacle and is thus unable to identify the direction of the steepest gradient.
In this case, the optimization process is stuck on the border of the domain.
In the third column of Fig.~\ref{fig:degeneracy}, i.e., case (c), rDSM corrects the degenerated simplex as introduced at the 9th iteration section, and ends at point (0.9952, -0.9860) with $J=0.0047$ after 100 iterations.
Comparing the above three cases, rDSM can identify the global optimum in the scenarios where DSM fails.

\begin{figure}[ht]
\centering
\includegraphics[width=0.7\textwidth]{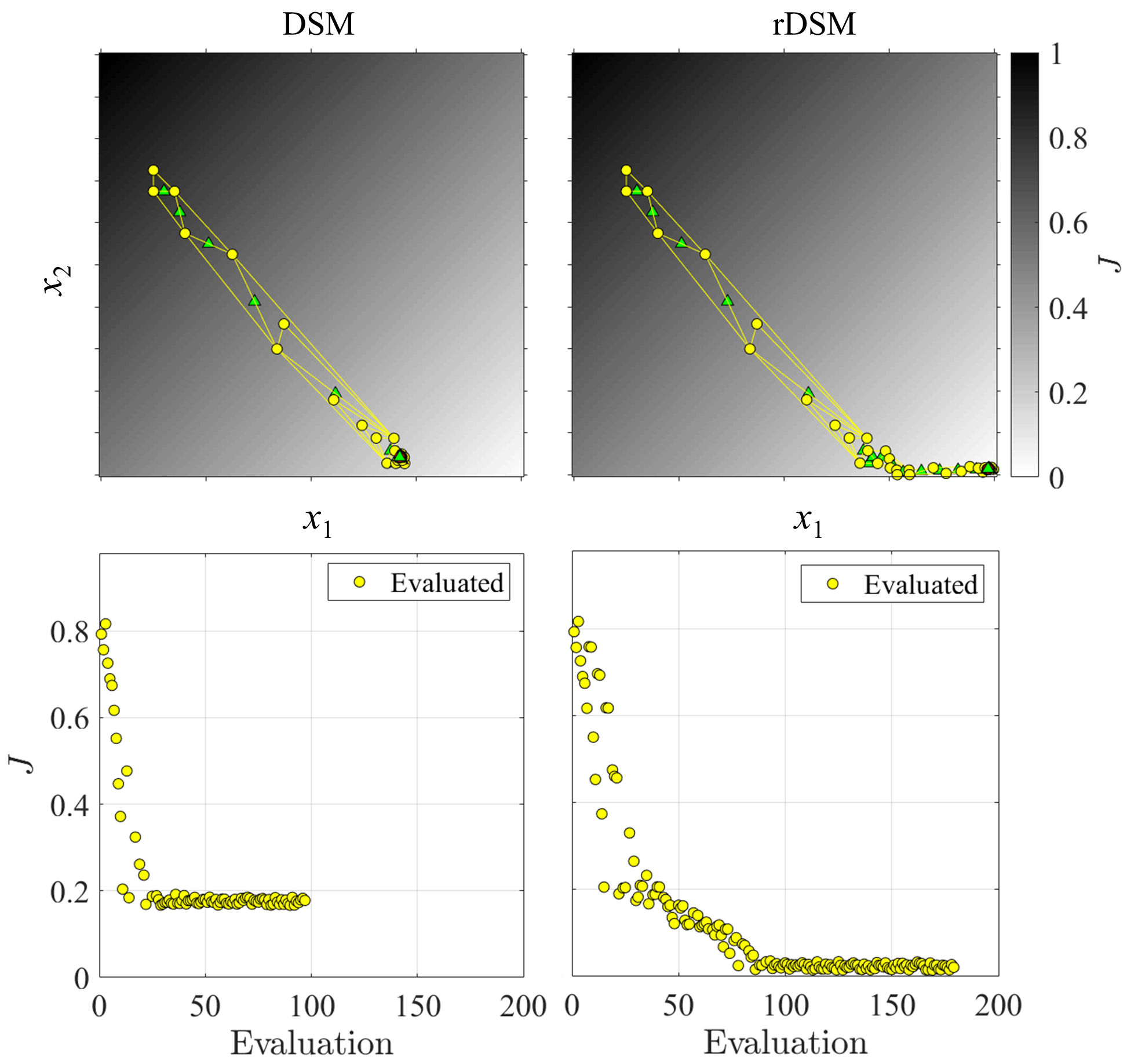}
\caption{Comparing the optimization process and corresponding learning curve of noisy linear gradient without obstacle by DSM and rDSM.
}
\label{fig:reevaluation}
\end{figure}

To demonstrate the robustness of rDSM to noise, we add a stochastic perturbation to the cost function.
The perturbation is bounded within the interval $[0,0.02]$~\cite{Nicholson2020, Krink2004,Taghiyeh2016}, and it is a uniform noise.
The number of iterations is limited to 50, maintaining the same starting point as in Fig.~\ref{fig:degeneracy}.
As illustrated in Fig.~\ref{fig:reevaluation}, the rDSM algorithm converges to a superior optimum compared to DSM.
Specifically, DSM ends at point (0.4216,-0.9148) with an objective value of 0.1780. 
In contrast, the optimum of rDSM is (0.9701,-0.9709) with $J=0.0212$.
This reevaluation mechanism facilitated rDSM in exploring the solution space more effectively, thereby identifying a more favorable optimum.
\begin{table}[h]
\centering
\caption{Comparing the results of minimizing the 2-dimensional problem by DSM and rDSM for repeating 20 times.}\label{table: Noise}
\begin{tabular}{cccc}
\toprule
Algorithm    & Type of noise  & DSM  & rDSM \\ \midrule
End point  & Uniform & \makecell{(0.9523 $\pm$0.0129, \\ -0.9530$\pm9\times10^{-4}$)}  &  \makecell{(0.9879$\pm 8\times 10^{-5}$,\\-0.9863$\pm3\times10^{-4}$) }                         \\
 $J $   & $\mathcal{U}[0,0.01]$ & 0.0248 $\pm 7\times10^{-4}$  &    0.0075  $\pm 2\times10^{-5}$                      \\ \hline
End point  & Uniform & \makecell{(0.5995$\pm 0.0252$,\\-0.9710$\pm 0.0016$)} &  \makecell{(0.9621$\pm 0.0057$,\\ -0.9629$\pm 8\times10^{-4}$) }                          \\
$J$  &$\mathcal{U}[0,0.02]$& 0.1111 $\pm$ 0.0014  & 0.0239 $\pm$ 0.0003                      \\ \hline
End point  & Gaussian &\makecell{(0.7727$\pm$0.0726, \\ -0.9435 $\pm$0.0033) }  &  \makecell{(0.9798$\pm 4 \times10^{-4}$,\\-0.9614$\pm$0.0013) }                      \\
$J$  & $\mathcal{N}(0,0.005)$ & 0.0662 $\pm$0.0034  & 0.0098 $\pm$0.0001                          \\
\hline
End point  & Gaussian & \makecell{(0.5585 $\pm$0.0836,\\-0.9812$\pm 4\times10^{-4}$)}  &    \makecell{(0.9559$\pm$0.0015,\\-0.9432$\pm$0.0026)}                         \\
$J $   & $\mathcal{N}(0,0.01)$ & 0.10555 $\pm$ 0.0051  &0.0133 $\pm$0.0001                         \\
\bottomrule
\end{tabular}
\end{table}

A comprehensive noise-robustness evaluation of rDSM is carried out by altering the cost function $J$ with two distinct noise models. 
Additive Gaussian perturbations $\mathcal{N} (0, \sigma^2)$ are imposed at two variance levels, while uniform disturbances $\mathcal{U} (a, b)$ are introduced at two interval widths.
Each noisy instance is optimized 20 times from the same starting point (-0.75, 0.35); all other settings remain fixed so that only the stochastic contamination varies.
Table~\ref{table: Noise} presents the resulting sample means and variances of the converged minimum and corresponding costs.
In contrast, the DSM method consistently reaches a local minimizer, while the rDSM method advances towards the neighborhood of the true optimum, thereby confirming its superior stability under noise.

To calibrate the reevaluation trigger, we compare the thresholds $c^{s_i} \geq 1.5n$ and $c^{s_i} \geq 2n$ by the uniform noise case $\mathcal{U}[0,0.02]$.
With the 1.5 coefficient, rDSM requires 4.4 seconds on average. However, an average of 8.3 seconds is needed when we set $c^{s_i} \geq 2n$.
More importantly, in the case that the coefficient is 2, rDSM converges to (0.6298$\pm 0.0163$, -0.9837$\pm 0.0027$) instead of the real global minimum (1, -1). 
In contrast, when the coefficient is set as 1.5, rDSM can converge to (0.9621$\pm 0.0057$, -0.9629$\pm 8\times10^{-4}$), which is very close to (1, -1).
These results indicate that the setting $c^{s_i} \geq 1.5n$ offers the best compromise between computational expense and solution accuracy.

\subsection{High-dimensional problem}

In this section, we minimize the high-dimensional Rosenbrock function~\cite{BOUVRY2000, Shang2006, KIRTAS2024} to test the proposed rDSM algorithm.
The Rosenbrock function is also referred to as the valley function. Its global minimum is located in a narrow, parabolic valley, which poses a challenge to converge~\cite{Picheny2013}. The mathematical function can be expressed as follows:
\begin{equation} \label{eq:Nd ronsenbrock}
    J_r(\bm{x}) = \sum_{i=1}^{n-1}[100(x_{i+1}-x_i^2)^2+(x_i-1)^2],
\end{equation}
where, $n$ refers to the dimension of $J_r(\bm{x})$. 
$J_r(\bm{x})$ reaches global minimum at $\bm{x}=(1,1,\cdots,1)$ with value of 0. 
First, we set $n=5, x_i \in [-5,10]$ for $i = 1, 2, \dots, 5$. 
A random initial point (-0.9598, -1.66907, -0.19862, -3.61086, -3.77915) is provided for both DSM and rDSM; the comparison of results is shown in the Table.~\ref{table: 5D rosen}.
For the convergence stability, the cost is divided by 10000 during the optimization process.
The number of iterations is limited to 500.
The edge and volume thresholds are set as $1\times10^{-5}$.
\begin{table}[h]
\centering
\caption{Comparing results of minimizing the 5-dimensional Rosenbrock function by DSM and rDSM.}\label{table: 5D rosen}
\begin{tabular}{ccc}
\toprule
Algorithm     & DSM  & rDSM \\ \midrule
Initial point   &\makecell{(-0.9598, -1.6691, -0.1986,\\ -3.6109, -3.7792)}  & \makecell{(-0.9598, -1.6691, -0.1986,\\ -3.6109, -3.7792)}   \\
End point  & \makecell{(0.9537, 0.9042, 0.8259,\\ 0.6944, 0.4746)}   &  \makecell{($\bm{1.0000,1.0000,1.0000,}$\\ $\bm{1.0000,1.0000}$)}                            \\
\makecell{Number of \\Iterations}&  500   & 500             \\ 
\makecell{Number of \\evaluations}& 785    & 882              \\ 
$J_r$ & $1.66\times10^{-5}$  & $\bm{2.92\times10^{-10}}$                           \\
\bottomrule
\end{tabular}
\end{table}

The algorithm DSM converges after 785 evaluations of the cost function at the point (0.9537, 0.9042, 0.8259, 0.6944, 0.4746) with an associated cost of $1.66\times10^{-5}$.
In contrast, rDSM locates the global minimum, achieving a cost of $J = 2.92\times10^{-10}$, though it evaluates $J_r$ 97 times more than DSM.
Running time measurements for this five-dimensional Rosenbrock problem reveal that DSM requires 3.4 seconds, whereas rDSM requires 64.6 seconds. 
Although rDSM requires a substantially higher computational cost, its capacity to reach the global minimum is clearly superior to that of DSM.
Consequently, on the five-dimensional Rosenbrock test function, rDSM consistently identifies the global minimum, whereas classical DSM does not.

\section{Impact}\label{Sec: impact}

High-dimensional optimization problems pose significant challenges, particularly in scenarios where the gradient of the objective function is either unavailable or computationally expensive to calculate.
The rDSM introduced in this study demonstrates the potential to identify the optimum without a computationally intensive process associated with gradient calculations.
Furthermore, the reevaluation functionality makes rDSM more robust in online optimization based on experiments.
This functionality ensures objective recalculation 
when a point remains in the simplex for over 1.5$ n$ iterations.
thereby refining solution accuracy. 
Given these features, the rDSM framework is anticipated to be effectively applicable to a broad range of high-dimensional optimization problems.

The first application of rDSM will be optimizing the fan array wind generator (FAWG)~\cite{Songqi2024AIAA}, including 100 individually controllable fans. The optimization goal is to achieve uniform flow characteristics. 
The input parameters are the magnitudes of fan switching, and the objective function is the square of the absolute error between the measured wind speed and the target wind speed.
With the help of rDSM, the flow field generated by FAWG will achieve uniform wind which can be used to test the performance of unmanned aerial vehicles (UAVs).
rDSM is expected to be applied for larger FAWG, demonstrated in~\cite{Yutong2025PoF}.

\section{Conclusions}\label{Sec: conclusion}

This paper introduces a software implementation of the robust Downhill Simplex Method (rDSM), which is an advanced optimization algorithm. 
The rDSM enhances the Downhill Simplex Method (DSM) through the incorporation of degeneracy correction and reevaluation strategies.
Degeneracy correction is achieved by maximizing the simplex volume while preserving its perimeter, thus preventing simplex collapse. 
The reevaluation mechanism improves the convergence robustness for optimization problems with noise. 
The objective function value of a vertex will be reevaluated if this vertex remains in the optimization process over 1.5$n$, where $n$ is the dimension of the optimization space.
These enhancements enable the rDSM to a better convergence compared to DSM, and may be beneficial for high-dimensional problems, even in noisy environments. Yet, they do not solve the scaling problem in high-dimensional space.
The effectiveness of rDSM is demonstrated using a minimization problem involving a linear gradient with an obstacle, where the rDSM successfully locates the global minimum, whereas the classic DSM fails.
This study highlights the robustness and efficiency of the rDSM, positioning it as a valuable advancement in optimization algorithms for complex problems. 
Future work may focus on exploring the potential applications of the rDSM in diverse optimization scenarios and developing dimensionality-reduction techniques or hybrid machine learning solvers to extend the algorithm to very high dimensions.

\section*{Acknowledgements}
\begin{sloppypar}
We are grateful for the inspiring discussion with Dr. Guy Y. Cornejo Maceda about formulating the idea of rDSM, exploring coding techniques, and offering insightful writing suggestions. 
Bernd Noack acknowledges support
by the National Science Foundation of China (NSFC)
through grant 12172109, 
by Guangdong province, China, 
via the Natural Science and Engineering grant 2022A1515011492,
and by the Shenzhen Science and Technology Program under grants JCYJ20220531095605012, KJZD20230923115210021 and 29853MKCJ202300205.
Xiaozhou He would like to acknowledge the support by the National Natural Science Foundation of China (No 12372216), the Science, Technology, and Innovation Commission of Shenzhen Municipality (Nos. GXWD20220818113020001 and JCYJ20240813104853070).
\end{sloppypar}

\section*{Data Availability Statement}
The data that support the findings of this study are available from the corresponding author upon reasonable request.

\bibliographystyle{unsrt}
\bibliography{aipsamp0}










\end{document}